\documentclass[reqno]{amsart}

\usepackage{amssymb}
\usepackage{hyperref}
\usepackage{amsmath}
\usepackage{amsfonts}
\usepackage{amssymb}
\usepackage{graphicx}

\begin{document}
\title[ Positive Solutions of Nonlinear Problems at Resonance ]{Some remarks on positive solutions of nonlinear problems at
resonance}
\author[F. Haddouchi, S. Benaicha]{Faouzi Haddouchi, Slimane Benaicha}
\address{Faouzi Haddouchi\\
Department of Physics, University of Sciences and Technology of
Oran, El Mnaouar, BP 1505, 31000 Oran, Algeria}
\email{haddouch@univ-usto.dz}
\address{Slimane Benaicha \\
Department of Mathematics, University of Oran, Es-senia, 31000 Oran,
Algeria} \email{slimanebenaicha@yahoo.fr} 

\begin{abstract}
The proof of a result of J. J. Nieto \cite{Niet2} appeared in
\textquotedblleft Acta Math, Hung.\textquotedblright(1992)
concerning the positive solutions of nonlinear problems at resonance
is corrected and improved.
\end{abstract}

\maketitle \numberwithin{equation}{section}
\newtheorem{theorem}{Theorem}[section]
\newtheorem{lemma}[theorem]{Lemma}
\newtheorem{proposition}[theorem]{Proposition}
\newtheorem{corollary}[theorem]{Corollary}
\newtheorem{remark}[theorem]{Remark}

\section{Introduction}
The Method of differential inequalities or the method of upper and
lower solutions has been used by Nieto \cite{Niet2} to show the
existence of positive periodic solutions for a second order
nonlinear differential equation. Nieto \cite{Niet2} has obtained two
existence results of positive and negative solutions for a class of
nonlinear problems at resonance. However we would like to point out
that the proof of the first main result (Theorem 6) in \cite{Niet2}
is not correct. We also  improve Theorem 7 of \cite{Niet2}. The
correction of the proof of Theorem 6 in \cite{Niet2} is the
motivation of this brief paper.

\section{Positive solutions and the method of upper and lower solutions}
J. J. Nieto in the paper \cite{Niet2} studied the existence of
positive periodic solutions of the equation
\begin{equation}  \label{eq-1}
{u''}+u+\mu u^{2}=h(t)\text{, }u(0)=u(\tau )\text{, }{u'}(0)%
={u'}(\tau )\text{,}
\end{equation}

where $h(t)=\epsilon \cos \omega t$ is $\tau =2\pi \omega ^{-1}$
periodic, $ \mu \neq 0,$ $\epsilon \neq 0$ and $\omega >0$.

Nieto and Rao in \cite{Niet:Rao} gave the following result:

\begin{theorem} \label{theo 2.1}
Equation \eqref{eq-1} has a periodic solution if  \ $ 4\left\vert
\mu \right\vert .\left\vert \epsilon \right\vert <1.$
\end{theorem}

Making $\ s=\omega t$, \eqref{eq-1} becomes%
\begin{equation}  \label{eq-2}
{u''}+\omega ^{-2}\left[ u+\mu u^{2}-\epsilon \cos s\right] =0%
\text{, }u(0)=u(2\pi )\text{, }{u'}(0)={u'}(2\pi )\text{,%
}
\end{equation}

where $u=u(s)$ \ and $\ {u''}=\frac{d^{2}u}{ds^{2}}$.

\ \

Thus we are interested in the existence of $2\pi $-periodic solutions of %
\eqref{eq-2} and note that it is of the form%
\begin{equation}  \label{eq-3}
\left\{
\begin{array}{l}
-{u ''}(t)=f(t,u),\text{ \newline
}t\in \left[ 0,2\pi \right] , \\
u(0)=u(2\pi )\text{, }{u'}(0)={u'}(2\pi ).%
\end{array}%
\right.
\end{equation}

As usual, we say that $\ \alpha \in C^{2}(\left[ 0,2\pi \right]
,\mathbb{R})$ is a lower solution of \eqref{eq-3} if
\begin{align*}
-{\alpha ''}(t) \leq f(t,\alpha (t))\ \text{for}\ t\in \left[ 0,2\pi
\right] ,
\end{align*}
\begin{align*}
\alpha (0)=\alpha (2\pi )\ \text{and}\ {\alpha '}(0)\geq {\alpha
'}(2\pi ).
\end{align*}

Similarly, $\ \beta \in C^{2}(\left[ 0,2\pi \right] ,\mathbb{R} )$
is an upper solution of \eqref{eq-3} if
\begin{align*}
\ -{\beta ''}(t)\geq %
f(t,\beta (t)) \ \text{for}\ t\in \left[ 0,2\pi \right] ,
\end{align*}
\begin{align*}
\beta (0)=\beta (2\pi )\ \text{and}\ {\beta '}(0)\leq {\beta '}(2\pi
).
\end{align*}

\begin{theorem}\label{theo 2.2} \cite{Kan}.
If \eqref{eq-3} has an upper solution $\beta $ and a lower
solution $\alpha $ such that $\alpha \leq \beta $ in $\left[ 0,2\pi \right] $%
, then there exists at least one solution $u$ of \eqref{eq-3} with
$\alpha \leq u\leq \beta $ in $\left[ 0,2\pi \right].$
\end{theorem}

We are now in a position to prove the following result due to Nieto
\cite{Niet2} and then we critically observe that it corrects the
proof of Theorem 6 of \cite{Niet2} and improves it since we do not
impose any condition on the sign of the real parameter $\epsilon $.

\begin{theorem}\label{theo 2.3} If $\mu <0$ and $\ 4\left\vert \mu \right\vert
.\left\vert \epsilon \right\vert <1$, then there exists a positive
($2\pi \omega ^{-1}$)-periodic solution of \eqref{eq-1}.

\end{theorem}

\begin{proof}
Note that equation \eqref{eq-2} can be written in the form%
\begin{equation}  \label{eq-4}
-{u ''}(s)=f(s,u),
\end{equation}

where $f(s,u)=\omega ^{-2}\left[ u+\mu u^{2}-\epsilon \cos
s\right]$.

For all arbitrary $\epsilon \neq 0$ and $\mu <0$, let
$0<a_{2}<a_{1}$ be the real roots of \ $\mu a^{2}+a-\left\vert
\epsilon \right\vert =0$ \ and \ $b_{2}<0<b_{1}$ the real roots of
$\mu b^{2}+b+\left\vert \epsilon \right\vert =0$. \ Note that $\
b_{2}<0<a_{2}<a_{1}<b_{1}$.

Choose $r\in \left[ a_{2},a_{1}\right] $ and $\ R\geq b_{1}$ and define $%
\alpha (s)=r$ and \ $\beta (s)=R$ \ ($r<R$) for $\ s\in \left[
0,2\pi \right] $. Since $\ -\left\vert \epsilon \right\vert \leq
-\epsilon \cos s\leq \left\vert \epsilon \right\vert $; we obtain

\begin{align*}
f(s,\beta (s))&=\omega ^{-2}(R+\mu R^{2}-\epsilon \cos s)\\
&\le \omega ^{-2}(R+\mu R^{2}+\left\vert \epsilon \right\vert )\\
&\le 0=-{%
\beta ''}(s),
\end{align*}

\begin{align*}
f(s,\alpha (s))&=\omega ^{-2}(r+\mu r^{2}-\epsilon \cos s)\\
&\ge \omega ^{-2}(r+\mu r^{2}-\left\vert \epsilon \right\vert )\\
&\ge 0=-{%
\alpha ''}(s).
\end{align*}
Therefore, by Theorem \ref{theo 2.2}, there exists a solution $u$ of
\eqref{eq-2} such that $u\geq r>0$. This complete the proof.

\end{proof}

Now, we shall improve Theorem 7 in \cite{Niet2} since we do not
require $\epsilon <0$.

\begin{theorem}\label{theo 2.4}
If $\mu >0$ and $\ 4\left\vert \mu \right\vert .\left\vert \epsilon
\right\vert <1$, then \eqref{eq-1} has a negative ($2\pi \omega
^{-1}$)-periodic solution.
\end{theorem}

\begin{proof}
The same argument as in Theorem 7 of \cite{Niet2} will be used.

\end{proof}


\begin{thebibliography}{99}

\bibitem{Kan} R. Kannan, and V. Lakshmikantham, Existence of periodic solutions of nonlinear
boundary value problems and the method of upper and lower solutions,
Appl. Anal., {\bf17} (1984), 103-113.


\bibitem{Niet:Rao} J. J. Nieto and V. S. H. Rao, Periodic solutions of
second order nonlinear differential equations, Acta Math. Hungar.,
{\bf48} (1986), 59-66.


\bibitem{Niet2} J. J. Nieto, Positive solutions of nonlinear problems at
resonance, Acta Math. Hungar., {\bf59} (3-4) (1992), 339-344.

\end{thebibliography}
\end{document}